\documentclass[10pt]{amsart}
\usepackage[leqno]{amsmath}
\usepackage{amssymb,latexsym,soul,cite,amsthm,color,enumitem,graphicx,mathtools,microtype,accents}
\usepackage[colorlinks=true,urlcolor=green(munsell),citecolor=green(munsell),linkcolor=green(munsell),linktocpage,pdfpagelabels,bookmarksnumbered,bookmarksopen]{hyperref}
\definecolor{green(munsell)}{rgb}{0.0, 0.66, 0.47}
\usepackage[english]{babel}
\usepackage[left=2.5cm,right=2.5cm,top=2.5cm,bottom=2.5cm]{geometry}

\numberwithin{equation}{section}

\newtheorem{theorem}{Theorem}[section]
\theoremstyle{plain}
\newtheorem{lemma}[theorem]{Lemma}
\theoremstyle{plain}
\newtheorem{proposition}[theorem]{Proposition}
\theoremstyle{plain}

\theoremstyle{definition}
\newtheorem{remark}[theorem]{Remark}
\newtheorem{example}[theorem]{Example}

\newcommand{\N}{{\mathbb N}}

\newcommand{\R}{{\mathbb R}}
\newcommand{\eps}{\varepsilon}
\newcommand{\beq}{\begin{equation}}
\newcommand{\eeq}{\end{equation}}
\renewcommand{\le}{\leqslant}
\renewcommand{\ge}{\geqslant}

\newcommand{\restr}[2]{\left.#1\right|_{#2}}
\newcommand{\w}{W^{s,p}_0(\Omega)}
\newcommand{\fpl}{(-\Delta)_p^s\,}
\newcommand{\ds}{{\rm d}_\Omega^s}
\newcommand{\cs}{C_s^0(\overline{\Omega})}

\makeatletter
\newcommand{\lenomode}{\tagsleft@true}
\newcommand{\reqnomode}{\tagsleft@false}
\makeatother

\newenvironment{enumroman}{\begin{enumerate}

}{\end{enumerate}}

\title[Fractional $p$-Laplacian with noncoercive energy]{Five solutions for the fractional $p$-Laplacian\\ with noncoercive energy}

\author[S.\ Frassu, A.\ Iannizzotto]{Silvia Frassu, Antonio Iannizzotto}

\address[S.\ Frassu, A.\ Iannizzotto]{Department of Mathematics and Computer Science
\newline\indent
University of Cagliari
\newline\indent
Via Ospedale 72, 09124 Cagliari, Italy}
\email{silvia.frassu@unica.it, antonio.iannizzotto@unica.it}


\begin{document}

\begin{abstract}
We deal with a Dirichlet problem driven by the degenerate fractional $p$-Laplacian and involving a nonlinear reaction which satisfies, among other hypotheses, a $(p-1)$-linear growth at infinity with non-resonance above the first eigenvalue. The energy functional governing the problem is thus noncoercive. Thus we focus on the behavior of the reaction near the origin, assuming that it has a $(p-1)$-sublinear growth at zero, vanishes at three points, and satisfies a reverse Ambrosetti-Rabinowitz condition. Under such assumptions, by means of critical point theory and Morse theory, and using suitably truncated reactions, we show the existence of five nontrivial solutions: two positive, two negative, and one nodal.
\end{abstract}

\maketitle

\begin{center}
Version of \today\
\end{center}

\section{Introduction}\label{sec1}

\noindent
Nonlinear boundary value problems driven by nonlocal operators of fractional order are the subject of a vast and rapidly developing literature, both for their intrinsic mathematical interest and for their numerous applications (see for instance \cite{BV,DPV}). Most authors have focused on the semilinear case, i.e, on problems driven by the fractional Laplacian or variants of it with several kernels (see \cite{RO} and the monograph \cite{MRS}).
\vskip2pt
\noindent
In the quasilinear case, which has the fractional $p$-Laplacian as its prototype operator, the fundamental theory is still developing (for an introduction to the subject see \cite{MS,P}). The simplest problem to be studied in such framework is the following Dirichlet problem:
\beq\label{dir}
\begin{cases}
\fpl u = f(x,u)& \text{in $\Omega$} \\
u=0 & \text{in $\Omega^c$.}
\end{cases}
\eeq
Here $\Omega\subset\R^N$ ($N\ge 2$) is a bounded domain with $C^{1,1}$ boundary, $p \ge 2$, $s\in(0,1)$ s.t.\ $N>ps$, and the leading operator is the degenerate fractional $p$-Laplacian, defined for all $u:\R^N\to\R$ smooth enough and all $x\in\R^N$ by
\beq \label{fpl}
\fpl u(x)=2\lim_{\eps\to 0^+}\int_{B_\eps^c(x)}\frac{|u(x)-u(y)|^{p-2}(u(x)-u(y))}{|x-y|^{N+ps}}\,dy
\eeq
(which in the linear case $p=2$ reduces to the fractional Laplacian, up to a dimensional constant). The reaction $f:\Omega\times\R\to\R$ is a Carath\'eodory mapping, subject to a subcritical growth condition on the real variable. The current literature on problem \eqref{dir} stems from the study of the corresponding nonlinear eigenvalue problem in \cite{BP,LL}, and it includes existence results based on Morse theory \cite{ILPS}, as well as regularity results \cite{BL,IMS,IMS1}, a version of the 'Sobolev vs.\ H\"older minima' result \cite{IMS2}, maximum and comparison principles \cite{DQ,J}, sub-supersolutions \cite{FI,KKP}, and existence/multiplicity results for several types of reactions based on both topological and variational methods \cite{A,BM,CMS,FI1,IL,IMP,PSY}.
\vskip2pt
\noindent
In dealing with nonlinear elliptic boundary value problems, the most delicate case is when the reaction is asymptotically $(p-1)$-linear at infinity, namely, when uniformly for a.e.\ $x\in\Omega$
\beq\label{lin}
-\infty < \liminf_{|t|\to\infty}\frac{f(x,t)}{|t|^{p-2}t} \le \limsup_{|t|\to\infty}\frac{f(x,t)}{|t|^{p-2}t} < \infty,
\eeq
in particular when such limits lie above the principal eigenvalue of the leading operator. Indeed, such case prevents both coercivity of the energy functional (as in the sublinear case), and the Ambrosetti-Rabinowitz condition at infinity (which is typically required in the superlinear case and allows for the use of the mountain pass theorem, see \cite{AR}). Usually one requires some type of non-resonance condition, i.e., that the limits in \eqref{lin} do not coincide with any eigenvalue (see \cite{AZ}). Under such additional assumption, the existence of nontrivial solutions can be proved either via topological or variational methods (see \cite{ILPS,FI1}, respectively, for problem \eqref{dir}). Several results have been proved for (local) $p$-Laplacian problems with noncoercive energy, see for instance \cite{FP,PR}.
\vskip2pt
\noindent
In this paper, we apply the state-of-the-art theory on the fractional $p$-Laplacian to a quite specific situation, taking inspiration from an interesting work of Papageorgiou-Smyrlis \cite{PS} dealing with the $p$-Laplacian. Our hypotheses are the following: $f(x,\cdot)$ is asymptotically $(p-1)$-linear at $\pm \infty$ with the asymptotic slopes in \eqref{lin} above the principal eigenvalue; also, $f(x,\cdot)$ is $(p-1)$-sublinear at $0$, crosses the $0$-line at $a_{-}<0<a_{+}$, and satisfies  a quasi-monotonicity condition. Under such hypotheses we prove that \eqref{dir} has at least four nontrivial constant sign solutions, two positive and two negative (Theorem \ref{css}). Adding a reverse  Ambrosetti-Rabinowitz condition near the origin, \eqref{dir} gains a fifth, nodal solution (Theorem \ref{nod}).
\vskip2pt
\noindent
Our approach is variational, based on critical point theory and Morse theory, and also makes a wide use of truncations, sub-supersolutions, and comparison principle (we will use the recent results of \cite{IMP}). Note that, while assuming that the limits in \eqref{lin} are above the principal eigenvalue, we do not assume non-resonance, thus turning our attention on the behavior of $f(x,\cdot)$ near the origin (in such sense we generalize the result of \cite{PS}). We confine ourselves to the degenerate case $p\ge 2$ essentially due to regularity reasons (see \cite{IMS1}), while the variational machinery also works, with minor adjustments, for the singular case $1<p<2$. Also, we remark that our result is new, even in the semilinear case $p=2$ (fractional Laplacian).
\vskip2pt
\noindent
The paper has the following structure: in Section \ref{sec2} we recall the functional-analytic framework and some well-known results about fractional $p$-Laplacian problems; in Section \ref{sec3} we show the existence of constant sign solutions (two positive and two negative); and  in Section \ref{sec4} we prove the existence of a nodal solution.
\vskip4pt
\noindent
{\bf Notation:} Throughout the paper, for any $A\subset\R^N$ we shall set $A^c=\R^N\setminus A$. For any two measurable functions $f,g:\Omega\to\R$, $f\le g$ in $\Omega$ will mean that $f(x)\le g(x)$ for a.e.\ $x\in\Omega$ (and similar expressions). The positive (resp., negative) part of $f$ is denoted $f^+$ (resp., $f^-$). If $X$ is an ordered Banach space, then $X_+$ will denote its non-negative order cone. For all $r\in[1,\infty]$, $\|\cdot\|_r$ denotes the standard norm of $L^r(\Omega)$ (or $L^r(\R^N)$, which will be clear from the context). Every function $u$ defined in $\Omega$ will be identified with its $0$-extension to $\R^N$. The constant functions of the type $u(x)=a$ are simply denoted by $a$. Moreover, $C$ will denote a positive constant (whose value may change case by case).

\section{Preliminaries}\label{sec2}

\noindent
In this section introduce the functional-analytic framework for problem \eqref{dir} and recall some useful results about the fractional $p$-Laplacian.
\vskip2pt
\noindent
First, for any Banach space $(X,\|\cdot\|)$ with topological dual $(X^*,\|\cdot\|_*)$ and all $u_0\in X$, $\rho>0$ we set
\[\overline{B}_{\rho}(u_0)=\big\{u \in X: \|u-u_0\|\le \rho\big\}.\]
Let $\Phi\in C^1(X)$ be a functional. For all $c\in\R$ we set
\[\Phi^c = \big\{u\in X:\,\Phi(u)\le c\big\}.\]
We say that $\Phi$ satisfies the Palais-Smale condition at level $c\in\R$, shortly $(PS)_c$, if every sequence $(u_n)$ in $X$, s.t.\ $\Phi(u_n)\to c$ and $\Phi'(u_n)\to 0$ in $X^*$, admits a (strongly) convergent subsequence. We say that $\Phi$ satisfies $(PS)$, if $\Phi$ satisfies $(PS)_c$ for any $c\in\R$.
\vskip2pt
\noindent
We denote by $K(\Phi)$ the set of all critical points of $\Phi$, and by $K_c(\Phi)$ the set of all $u\in K(\Phi)$ s.t.\ $\Phi(u)=c$. We say that $u\in K_c(\Phi)$ is an isolated critical point, if there exists a neighborhood $U\subset X$ of $u$ s.t.\ $K(\Phi)\cap U=\{u\}$, and in such case, for all $k\in\N$ we define the $k$-th critical group of $\Phi$ at $u$ as
\[C_k(\Phi,u) = H_k(\Phi^c\cap U,\Phi^c\cap U\setminus\{u\}),\]
where $H_k(\cdot,\cdot)$ denotes the $k$-th singular homology group for a topological pair (see \cite[Chapter 6]{MMP}).
\vskip2pt
\noindent
Following \cite{DPV}, for all measurable $u:\R^N\to\R$ we set
\[[u]_{s,p}^p=\iint_{\R^N\times\R^N}\frac{|u(x)-u(y)|^p}{|x-y|^{N+ps}}\,dx\,dy,\]
and define the fractional Sobolev spaces
\[W^{s,p}(\R^N)=\big\{u\in L^p(\R^N):\,[u]_{s,p}<\infty\big\},\]
\[\w=\big\{u\in W^{s,p}(\R^N):\,u(x)=0 \ \text{in $\Omega^c$}\big\},\]
the latter being a uniformly convex, separable Banach space with norm $\|u\|=[u]_{s,p}$. The dual space of $\w$ is denoted $W^{-s,p'}(\Omega)$, with norm $\|\cdot\|_{-s,p'}$. The embedding $\w\hookrightarrow L^q(\Omega)$ is continuous for all $q\in[1,p^*_s]$ and compact for all $q\in[1,p^*_s)$, with  $p_s^*= Np / (N-ps)$.
\vskip2pt
\noindent
Also we recall from \cite[Definition 2.1]{IMS} the following space
\[\widetilde{W}^{s,p}(\Omega)=\Big\{u\in L_{\rm loc}^p(\R^N):\, \exists \, \Omega' \Supset  \Omega \text{ s.t. } u \in W^{s,p}(\Omega'), \int_{\R^N} \frac{|u(x)|^{p-1}}{(1+|x|)^{N+ps}}\,dx < \infty\Big\}.\]
Clearly, any constant function $u(x)=a$ in $\R^N$ lies in $\widetilde{W}^{s,p}(\Omega)$. Further, $\w\subseteq\widetilde{W}^{s,p}(\Omega)$, and vice versa we have:

\begin{lemma}\label{zero}
If $u \in \widetilde{W}^{s,p}(\Omega)$ and $u=0$ in $\Omega^c$, then $u \in \w$.
\end{lemma}
\begin{proof}
By assumption, $u\in L^p(\R^N)$ and $u=0$ in $\Omega^c$, so there remains to show that $[u]_{s,p}<\infty$. Indeed, by symmetry we have
\begin{align*}
\iint_{\R^N\times\R^N}\frac{|u(x)-u(y)|^p}{|x-y|^{N+ps}}\,dx\,dy &= \iint_{\Omega'\times\Omega'}\frac{|u(x)-u(y)|^p}{|x-y|^{N+ps}}\,dx\,dy + 2\iint_{\Omega'\times\Omega'^c}\frac{|u(x)-u(y)|^p}{|x-y|^{N+ps}}\,dx\,dy\\
&+ \iint_{\Omega'^c\times\Omega'^c}\frac{|u(x)-u(y)|^p}{|x-y|^{N+ps}}\,dx\,dy.
\end{align*}
We note that the first integral on the right-hand side is finite since $u \in W^{s,p}(\Omega')$ and the third integral is zero by the condition on $\Omega^c$. So we focus on the second integral:
\begin{align*}
\iint_{\Omega'\times\Omega'^c}\frac{|u(x)-u(y)|^p}{|x-y|^{N+ps}}\,dx\,dy &=  \iint_{\Omega \times\Omega'^c}\frac{|u(x)|^p}{|x-y|^{N+ps}}\,dx\,dy\\
&\le \iint_{\Omega\times\Omega'^c}\frac{|u(x)|^p}{C(1+|y|^{N+ps})}\,dx\,dy \\
&=  \frac{1}{C}\int_{\Omega} |u(x)|^p\,dx \int_{\Omega'^c}\frac{1}{(1+|y|^{N+ps})}\,dy < \infty.
\end{align*}
Hence $[u]_{s,p}<\infty$, so $u \in \w$.
\end{proof}

\noindent
By \cite[Lemma 2.3]{IMS}, we can rephrase the fractional $p$-Laplacian as an operator $\fpl:  \widetilde{W}^{s,p}(\Omega) \to W^{-s,p'}(\Omega)$ defined for all 
$u \in \widetilde{W}^{s,p}(\Omega)$, $v \in \w$ by 
\[\langle \fpl u, v\rangle = \iint_{\R^N \times \R^N} \frac{|u(x)-u(y)|^{p-2} (u(x)-u(y)) (v(x)-v(y))}{|x-y|^{N+ps}}\,dx\,dy.\]
This definition is equivalent to \eqref{fpl}, provided $u$ is sufficiently smooth. The restricted operator $\fpl:  \w \to W^{-s,p'}(\Omega)$ is continuous, and satisfies the $(S)_+$-property, i.e., whenever $(u_n)$ is a sequence in $\w$ s.t.\ $u_n \rightharpoonup u$ in $\w$ and
\[\limsup_n \, \langle \fpl u_n, u_n-u \rangle \le 0,\]
then $u_n \to u$ in $\w$ (see \cite[Lemma 2.1]{FI}, \cite[Lemma 3.2]{FRS}). By \cite[Lemma 2.1]{IL}, for all $u\in\w$ we have
\beq \label{pnp}
\|u^{\pm}\|^p \le \langle \fpl u, \pm u^{\pm}\rangle.
\eeq
Moreover, the operator $\fpl$ is strictly $(T)$-monotone, namely:

\begin{proposition}\label{stm}
{\rm \cite[proof of Lemma 9]{LL}} Let $u, v \in \widetilde{W}^{s,p}(\Omega)$ s.t.\ $(u-v)^+ \in \w$ satisfy 
\[\langle \fpl u - \fpl v, (u-v)^+\rangle \le 0,\]
then $u \le v$ in $\Omega$. 
\end{proposition}

\noindent
We consider problem \eqref{dir} under the following hypothesis:
\begin{itemize}[leftmargin=1cm]
\item[${\bf H}_0$] $f:\Omega\times\R\to\R$ is a Carath\'{e}odory function, and there exist $c_0>0$, $r\in(1,p^*_s)$ s.t.\ for a.e.\ $x\in\Omega$ and all $t\in\R$
\[|f(x,t)| \le c_0 (1+|t|^{r-1}).\]
\end{itemize}
We say that $u\in \widetilde{W}^{s,p}(\Omega)$ is a {\em (weak) supersolution} of \eqref{dir}, if $u \ge 0$ in $\Omega^c$ and for all $v\in\w_+$
\[\langle \fpl u, v\rangle \ge \int_\Omega f(x,u)v\,dx.\]
Similarly, we define a {\em (weak) subsolution} of \eqref{dir}. We say that $(\underline{u}, \overline{u}) \in \widetilde{W}^{s,p}(\Omega) \times \widetilde{W}^{s,p}(\Omega)$
is a {\em sub-supersolution pair} of \eqref{dir}, if $\underline{u}$ is a subsolution, $\overline{u}$ is a supersolution, and $\underline{u} \le \overline{u}$ in $\Omega$. Finally, $u \in \w$ is a {\em (weak) solution} of \eqref{dir}  if $u$ is both a super- and a subsolution, i.e., if for all $v \in \w$
\[\langle \fpl u, v\rangle = \int_\Omega f(x,u)v\,dx.\]
On spaces $\w, \widetilde{W}^{s,p}(\Omega)$ we consider the pointwise partial ordering, inducing a lattice structure. Set
\[u \land v = \min\{u,v\}, \ u \vee v =\max\{u,v\}.\]
We recall that  the pointwise minimum of supersolutions is a supersolution, as well as the maximum of subsolutions is a subsolution (see also \cite{KKP}):

\begin{proposition}\label{sss}
{\rm\cite[Lemma 3.1]{FI}} Let ${\bf H}_0$ hold and $u_1,u_2 \in \widetilde{W}^{s,p}(\Omega)$:
\begin{enumroman}
\item \label{sss1} if $u_1, u_2$ are supersolutions of \eqref{dir}, then so is $u_1 \land u_2$;
\item \label{sss2} if $u_1, u_2$ are subsolutions of \eqref{dir} then so is $u_1 \vee u_2$.
\end{enumroman}
\end{proposition}

\noindent
For any sub-supersolution pair $(\underline{u}, \overline{u})$ of \eqref{dir} we define the solution set
\[\mathcal{S}(\underline{u},\overline{u}) = \big\{u \in \w:  u \text{ solves } \eqref{dir}, \, \underline{u} \le u \le \overline{u} \text{ in } \Omega\big\}.\]
The properties of $\mathcal{S}(\underline{u},\overline{u})$ were studied in \cite{FI}, and amount at the following:

\begin{proposition}\label{ext}
{\rm\cite[Theorem 3.5]{FI}} Let ${\bf H}_0$ hold, $(\underline{u}, \overline{u})$ be a sub-supersolution pair of \eqref{dir}. Then $\mathcal{S}(\underline{u},\overline{u})$ is both upward and downward directed, compact in $\w$, and it contains a smallest and a biggest element.
\end{proposition}

\noindent
Some words about regularity of solutions are now in order. First, we have a uniform {\em a priori} bound:

\begin{proposition}\label{apb}
{\rm\cite[Theorem 3.3]{CMS}} Let ${\bf H}_0$ hold, $u \in \w$ be a solution of \eqref{dir}. Then, $u \in L^{\infty}(\Omega)$ with $\|u\|_{\infty} \le C$, for some $C=C(\|u\|)>0$.
\end{proposition}

\noindent
Solutions of fractional order equations have a good interior regularity, but they may fail to be smooth up to the boundary. Thus, an important role in regularity theory for nonlinear, nonlocal operators is played by the following weighted H\"older spaces, with weight
\[\ds(x)= \mathrm{dist}(x, \Omega^c)^s.\]
Set
\[\cs = \Big\{u \in C^0(\overline{\Omega}): \frac{u}{\ds} \ \text{has a continuous extension to} \ \overline{\Omega} \Big\},\]
and for all $\alpha \in (0,1)$
\[C_s^{\alpha}(\overline{\Omega})= \Big\{u \in C^0(\overline{\Omega}): \frac{u}{\ds} \ \text{has a $\alpha$-H\"older continuous extension to} \ \overline{\Omega}\Big\},\]
with norms, respectively,
\[\|u\|_{0,s}= \Big\|\frac{u}{\ds}\Big\|_{\infty}, \ \|u\|_{\alpha,s}= \|u\|_{0,s} + \sup_{x \neq y} \frac{|u(x)/\ds(x) - u(y)/\ds(y)|}{|x-y|^{\alpha}}.\]
The embedding $C_s^{\alpha}(\overline{\Omega}) \hookrightarrow \cs$ is compact for all $\alpha \in (0,1)$. By \cite[Lemma 5.1]{ILPS}, the positive cone $\cs_+$ of $\cs$ has a nonempty interior given by
\[{\rm int}(\cs_+)= \Big\{u \in \cs:\, \inf_{\Omega}\frac{u}{\ds} > 0\Big\}.\]
Combining Proposition \ref{apb} and \cite[Theorem 1.1]{IMS1}, we have the following global regularity result for the degenerate case $p\ge 2$:

\begin{proposition}\label{reg}
Let ${\bf H}_0$ hold, $u \in \w$ be a solution of \eqref{dir}. Then, $u \in C_s^{\alpha}(\overline{\Omega})$  for some $\alpha \in (0,s]$.
\end{proposition}

\noindent
We recall from \cite{IMP} two general maximum/comparison principles:

\begin{proposition}\label{max}
{\rm\cite[Theorem 2.6]{IMP}} Let $g \in C^0(\R) \cap BV_{\rm loc}(\R)$, $u \in \widetilde{W}^{s,p}(\Omega) \cap C^0(\bar{\Omega})\setminus\{0\}$ s.t.\ 
\[\begin{cases}
\fpl u +g(u) \ge g(0) & \text{in $\Omega$} \\
u\ge0 & \text{in $\R^N$.}
\end{cases}\]
Then, 
\[\inf_{\Omega} \frac{u}{\ds}>0.\]
In particular, if $u \in \cs$, then $u \in {\rm int}(\cs_+)$. 
\end{proposition}

\begin{proposition}\label{comp}
{\rm\cite[Theorem 2.7]{IMP}} Let $g \in C^0(\R) \cap BV_{\rm loc}(\R)$, $u \in \w \cap C^0(\overline\Omega)$, $v \in \widetilde{W}^{s,p}(\Omega) \cap C^0(\bar{\Omega})$, $u \neq v$ and $C>0$ s.t.\
\[\begin{cases}
\fpl u +g(u) \le \fpl v +g(v) \le C & \text{in $\Omega$} \\
0 < u\le v & \text{in $\Omega$} \\
v \ge 0 & \text{in $\Omega^c$.}
\end{cases}\]
Then, 
\[\inf_{\Omega} \frac{v-u}{\ds}>0.\]
In particular, if $u,v\in{\rm int}(\cs_+)$, then $v-u\in{\rm int}(\cs_+)$.
\end{proposition}

\noindent
Now we define an  energy functional for problem \eqref{dir} by setting for all $(x,t)\in\Omega\times\R$
\[F(x,t)=\int_0^t f(x,\tau)\,d\tau,\]
and for all $u\in \w$
\[\Phi(u)= \frac{\|u\|^p}{p} -  \int_{\Omega} F(x,u)\,dx.\]
By ${\bf H}_0$, it is easily seen that $\Phi \in C^1(\w)$ with G\^ateaux derivative given for all $u,v\in\w$ by
\[\langle\Phi'(u),v\rangle = \langle\fpl u,v\rangle-\int_\Omega f(x,u)v\,dx.\]
Clearly, then, critical points of $\Phi$ coincide with the solutions of \eqref{dir}. Also, $\Phi$ is sequentially weakly lower semicontinuous, and every bounded $(PS)$-sequence for $\Phi$ has a convergent subsequence (see \cite[Proposition 2.1]{ILPS}). Another useful property is that local minima of $\Phi$ in the topologies of $\w$ and $\cs$, respectively, coincide:

\begin{proposition}\label{svh}
{\rm \cite[Theorem 1.1]{IMS1}} Let ${\bf H}_0$ hold, $u\in\w$. Then, the following are equivalent:
\begin{enumroman}
\item\label{svh1} there exists $\rho>0$ s.t.\ $\Phi(u+v)\ge\Phi(u)$ for all $v\in\w\cap C_s^0(\overline\Omega)$, $\|v\|_{0,s}\le\rho$;
\item\label{svh2} there exists $\sigma>0$ s.t.\ $\Phi(u+v)\ge\Phi(u)$ for all $v\in\w$, $\|v\| \le\sigma$.
\end{enumroman}
\end{proposition}

\noindent
Finally, we recall some properties of the following weighted eigenvalue problem, with weight $m \in L^{\infty}(\Omega)_+ \setminus \{0\}$:
\beq\label{evm}
\begin{cases}
\fpl u = \lambda m(x) |u|^{p-2}u & \text{in $\Omega$} \\
u=0 & \text{in $\Omega^c$.}
\end{cases}
\eeq
Problem \eqref{evm} admits an unbounded, nondecreasing sequence of positive variational (Lusternik-Schnirelmann) eigenvalues $(\lambda_k(m))$. In particular, as a special case of \cite[Theorem 1.1]{I} (see also \cite[Proposition 3.4]{FRS}) we have:

\begin{proposition}\label{fev}
Let $m\in L^\infty(\Omega)_+\setminus\{0\}$. Then, the smallest eigenvalue of \eqref{evm} is
\[\lambda_1(m)= \inf_{u \in \w \setminus \{0\}} \frac{\|u\|^p}{\int_{\Omega} m(x)|u|^p\,dx}>0.\]
Besides,
\begin{enumroman}
\item\label{fev1} $\lambda_1(m)$ is simple, isolated, with constant sign eigenfunctions, while for any eigenvalue $\lambda > \lambda_1(m)$ of \eqref{evm} $\lambda$-eigenfunctions are nodal;
\item\label{fev2} if $\tilde{m} \in L^\infty(\Omega)_+\setminus\{0\}$ is s.t.\ $m\le \tilde{m}$ in $\Omega$ and $m\not\equiv \tilde{m}$, we have $\lambda_1(m)>\lambda_1(\tilde{m})$.
\end{enumroman}
When $m=1$, we set $\lambda_1(1)=\lambda_1$ and we denote by $e_1 \in {\rm int}(\cs_+)$ the unique, positive $\lambda_1$-eigenfunction s.t.\ $\|e_1\|_p=1$.
\end{proposition}

\section{Constant sign solutions}\label{sec3}

\noindent
In this section we study problem \eqref{dir} under the following hypotheses:

\begin{itemize}[leftmargin=1cm]
\item[${\bf H}_1$] $f:\Omega\times\R\to\R$ is a Carath\'{e}odory mapping satisfying
\begin{enumroman}
\item\label{h11} there exists $c_0>0$ s.t.\ for a.e.\ $x\in\Omega$ and all $t \in \R$
\[|f(x,t)| \le c_0(1+|t|^{p-1});\]
\item\label{h12} there exist $\eta_1, \eta_2 \in L^{\infty}(\Omega)$ s.t.\ $\eta_1 \ge \lambda_1$ in $\Omega$, $\eta_1 \not\equiv \lambda_1$, and uniformly for a.e.\ 
$x \in \Omega$
\[\eta_1(x) \le \liminf_{|t|\to \infty}\frac{f(x,t)}{|t|^{p-2}t} \le \limsup_{|t|\to \infty}\frac{f(x,t)}{|t|^{p-2}t} \le \eta_2(x);\] 
\item\label{h13} there exist $\delta_0, c_1>0$, $q\in (1,p)$ s.t.\ for a.e.\ $x \in \Omega$
\[f(x,t) \ \begin{cases}
\ge c_1 t^{q-1} & \text{for all $t \in [0,\delta_0]$,} \\
\le c_1 |t|^{q-2}t & \text{for all $t \in [-\delta_0,0]$;}
\end{cases}\] 
\item\label{h14} there exist $a_{-} ,a_{+}$ with $\min\{a_{+},-a_{-}\}>\delta_0$ s.t.\  for a.e.\ $x \in \Omega$
\[f(x,a_{-}) = f(x,0) = f(x,a_{+}) = 0;\]
\item\label{h15} there exists $c_2>0$ s.t.\ for a.e.\ $x\in\Omega$ the map
\[t \mapsto f(x,t)+c_2|t|^{p-2}t\]
is nondecreasing in $[a_{-},a_{+}]$.
\end{enumroman}
\end{itemize}

\noindent
Hypotheses ${\bf H}_1$ include a subcritical growth condition \ref{h11}, $(p-1)$-linear asymptotic growth with non-resonance on the first eigenvalue \ref{h12}, a $(p-1)$-sublinear behavior near the origin \ref{h13}, vanishing at $a_-<0<a_+$ (not necessarily with a sign change) \ref{h14}, and a quasi-monotonicity condition \ref{h15}. We note that, contrary to the assumptions of \cite{PS}, we do not require that $\eta_2<\lambda_2$ in $\Omega$, that is, we allow resonance on higher eigenvalues.
\vskip2pt
\noindent
Clearly, ${\bf H}_1$ imply ${\bf H}_0$ (with $r=p$), so all results of Section \ref{sec2} apply here, with $f$, $\Phi$ defined as above. Also, by ${\bf H}_1$ \ref{h14} problem \eqref{dir} admits the trivial solution $u=0$. Without loss of generality, we may assume that \eqref{dir} has only {\em finitely many} solutions.

\begin{example}\label{exa}
We present here an autonomous reaction $f\in C^0(\R)$ satisfying ${\bf H}_1$. Set $p=2$, and fix real numbers $\eta>\lambda_1$, $\gamma>\eta+1$, and set for all $t\ge 0$
\[f(t) = \frac{\eta t^2-\gamma t+\sqrt{t}}{t+1}.\]
We focus on the positive semiaxis. Clearly, $f$ satisfies \ref{h11}. Since
\[\lim_{t\to\infty}\frac{f(t)}{t} = \eta,\]
then we have \ref{h12}. Taking $q=3/2$, we have
\[\lim_{t\to 0^+}\frac{f(t)}{\sqrt{t}} = 1,\]
which implies \ref{h13}. Since $f(t)>0$ for all $t>0$ small enough and
\[f(1) = \frac{\eta-\gamma+1}{2} < 0,\]
then there exists $a_+\in(0,1)$ s.t.\ $f(a_+)=0$, whence \ref{h14}. Finally, since $f'(t)$ is bounded below in $(0,\infty)$, we can find $c_2>0$ s.t.\
\[t\mapsto f(t)+c_2t\]
is increasing \ref{h15}. Taking $f(t)=-f(-t)$ for all $t<0$, we complete the definition.
\end{example}

\noindent
We begin our study by proving the existence of a first positive solution:

\begin{lemma}\label{min}
If ${\bf H}_1$ hold, then \eqref{dir} has a solution $u_0 \in {\rm int}(\cs_+)$, which is a local minimizer of $\Phi$.
\end{lemma}
\begin{proof}
Set for all $(x,t) \in \Omega \times \R$
\[f_0(x,t) =
\begin{cases}
0 & \text{if $t<0$,} \\
f(x,t) & \text{if $0\le t \le a_+$,}\\
0 & \text{if $t>a_+$,}
\end{cases}\]
with $a_+>0$ as in ${\bf H}_1$ \ref{h14}. By ${\bf H}_1$ \ref{h11} \ref{h14} we see that $f_0: \Omega \times \R \to \R$ satisfies ${\bf H}_0$. Also, it is bounded. Accordingly, set for all $(x,t) \in \Omega \times \R$
\[F_0(x,t) = \int_0^t f_0(x,\tau)\,d\tau,\]
and for all $u\in\w$
\[\Phi_0(u)= \frac{\|u\|^p}{p} -  \int_{\Omega} F_0(x,u)\,dx.\]
As seen in Section \ref{sec2}, $\Phi_0 \in C^1(\w)$ is sequentially weakly lower semicontinuous. Besides, $\Phi_0$ is coercive. Indeed, for all $u\in\w$
\begin{align*}
\Phi_0(u) &\ge \frac{\|u\|^p}{p} -  \int_{\Omega} C |u|\,dx\\
&\ge \frac{\|u\|^p}{p} - C \|u\|,
\end{align*}
and the latter tends to $\infty$ as $\|u\| \to \infty$. So there exists $u_0 \in \w$ s.t.\ 
\beq\label{min1}
\Phi_0(u_0) = \inf_{u\in\w}\Phi_0(u) =: m_0.
\eeq
Let $e_1\in{\rm int}(C^0_s(\overline\Omega)_+)$ be defined as in Proposition \ref{fev}. In particular, $e_1\in L^\infty(\Omega)$, so for all $\tau >0$ small enough (recalling ${\bf H}_1$ \ref{h14}) we have in $\Omega$
\[0<\tau e_1\le\delta_0 < a_{+}.\]
Hence, by ${\bf H}_1$ \ref{h13} we have 
\begin{align*}
\Phi_0(\tau e_1) &\le \frac{\tau^p}{p}\|e_1\|^p-\int_\Omega\frac{c_1}{q}(\tau e_1)^q\,dx\\ 
&= \frac{\lambda_1}{p}\tau^p- \frac{c_1\|e_1\|_q^q}{q} \tau^q,
\end{align*}
and the latter is negative for all $\tau >0$ small enough (as $q<p$). So in \eqref{min1} we have
\[m_0 = \Phi_0(u_0) < 0 = \Phi_0(0),\]
hence $u_0 \neq 0$. By \eqref{min1} we have $\Phi'_0(u_0)=0$, i.e. weakly in $\Omega$
\beq\label{min2}
\fpl u_0 = f_0(x,u_0).
\eeq
Testing \eqref{min2} with $-u_0^- \in \w$ and using \eqref{pnp}, we get 
\begin{align*}
\|u_0^-\|^p &\le \langle \fpl u_0, -u_0^-\rangle \\
&= \int_\Omega f_0(x,u_0) (-u_0^-)\,dx = 0,
\end{align*}
so $u_0 \ge 0$ in $\Omega$. Besides, note that $a_+ \in \widetilde{W}^{s,p}(\Omega)$ satisfies weakly in $\Omega$
\[\fpl a_+ =0.\]
By Lemma \ref{zero} we have $(u_0-a_+)^+ \in \w$. By using such function as test in \eqref{min2} and the above equation, we get
\[\langle \fpl u_0 - \fpl a_+, (u_0-a_+)^+\rangle = \int_\Omega f_0(x,u_0) (u_0-a_+)^+\,dx =0.\] 
By Proposition \ref{stm} we have $0\le u_0 \le a_+$ in $\Omega$. Thus, in \eqref{min2} we can replace $f_0$ with $f$ and see that $u_0\in\w\setminus\{0\}$ is a solution of \eqref{dir}. By Proposition \ref{reg} (with $g(t)=c_2|t|^{p-2}t$) we have $u_0 \in C_s^{\alpha}(\overline{\Omega})\setminus\{0\}$. Further, by ${\bf H}_1$ \ref{h15} we have weakly in $\Omega$ 
\[\begin{cases}
\fpl u_0 + c_2 u_0^{p-1} = f(x,u_0) + c_2 u_0^{p-1} \ge f(x,0)=0 \\
u_0 \ge 0. 
\end{cases}\]
By Proposition \ref{max} (with the same choice of $g$), we deduce $u_0 \in {\rm int}(\cs_+)$.
\vskip2pt
\noindent
Finally, we prove that $u_0$ is a local minimizer of $\Phi$. By ${\bf H}_1$ \ref{h14} \ref{h15} we have weakly in $\Omega$
\[\begin{cases}
\fpl u_0 + c_2 u_0^{p-1} \le  f(x,a_+) + c_2 a_+^{p-1}= \fpl a_+ + c_2 a_+^{p-1} = c_2 a_+^{p-1}\\
0 < u_0 \le a_+, 
\end{cases}\]
while $u_0\in\w$ clearly implies $u_0\neq a_+$. So, by Proposition \ref{comp} (still with the same $g$), we infer that 
\[\inf_{\Omega} \frac{a_+-u_0}{\ds}>0,\]
in particular $0<u_0 < a_+$ in $\Omega$. Set
\[\mathcal{U}= \big\{u \in \w \cap \cs: 0 < u < a_+ \ \text{in $\Omega$}\big\},\]
an open set in $\cs$ s.t.\ $u_0 \in \mathcal{U}$. By \eqref{min1} we have for all $u \in \mathcal{U}$
\[\Phi(u)=\Phi_0(u) \ge \Phi_0(u_0)=\Phi(u_0).\]
So, $u_0$ is a local minimizer of $\Phi$ in $\cs$. By Proposition \ref{svh}, $u_0$ is a local minimizer of $\Phi$ in $\w$ as well.
\end{proof}

\noindent
Exploiting the asymptotic behavior of $f(x,\cdot)$ at $\infty$, we find a second positive solution:

\begin{lemma}\label{sps}
If ${\bf H}_1$ hold, then \eqref{dir} has a solution $u_1 \in {\rm int}(\cs_+)$ s.t.\ $u_1 \ge u_0$, $u_1 \neq u_0$.
\end{lemma}
\begin{proof}
Set for all $(x,t) \in \Omega \times \R$
\[f_1(x,t) =
\begin{cases}
f(x,u_0) & \text{if $t<u_0(x)$,} \\
f(x,t) & \text{if $t \ge u_0(x)$,}
\end{cases}\]
with $u_0 \in {\rm int}(\cs_+)$ given by Lemma \ref{min}. By ${\bf H}_1$ \ref{h11}, $f_1:\Omega\times\R\to\R$ satisfies ${\bf H}_0$. So, for all $(x,t)\in\Omega\times\R$ we set
\[F_1= \int_0^t f_1(x,\tau)\,d\tau.\]
We define $\Phi_1\in C^1(\w)$ by setting for all $u\in\w$
\[\Phi_1(u)= \frac{\|u\|^p}{p} - \int_{\Omega} F_1(x,u)\,dx.\]
We prove first that $\Phi_1$ satisfies $(PS)$. Indeed, let $(w_n)$ be a sequence in $\w$ s.t.\ $|\Phi_1(w_n)| \le C$ and $\Phi'_1(w_n) \to 0$ in $W^{-s,p'}(\Omega)$, i.e., there exists a sequence $(\eps_n)$ in $\R_0^+$ s.t.\ $\eps_n \to 0^+$ and for all $n \in \N$, $v \in \w$
\beq \label{sps1}
\Big |\langle \fpl w_n, v \rangle - \int_\Omega f_1(x,w_n)v\,dx \Big| \le \eps_n \|v\|.
\eeq
We aim at proving that $(w_n)$ is bounded in $\w$. Arguing by contradiction, assume that up to a subsequence $\|w_n\|\to\infty$. Testing \eqref{sps1} with $-w_n^- \in \w$ and applying \eqref{pnp}, we have for all $n\in\N$
\begin{align*}
\|w_n^-\|^p &\le \langle \fpl w_n, -w_n^- \rangle \\
&\le \int_\Omega f_1(x,w_n)(-w_n^-)\,dx + \eps_n \|w_n^-\| \\
&= \int_\Omega f(x,u_0)(-w_n^-)\,dx + \eps_n \|w_n^-\|\\
&\le C \|w_n^-\|,
\end{align*}
where we have used the definition of $f_1$ and ${\bf H}_1$ \ref{h11}. Since $p>1$, $(w_n^-)$ is bounded in $\w$. So
\[\|w_n^+\|= \|w_n+w_n^-\| \ge \|w_n\| - \|w_n^-\| \to \infty.\]
Set for all $n\in\N$
\[\hat{w}_n=\frac{w_n}{\|w_n\|},\]
so $\hat{w}_n\in\w$ with $\|\hat w_n\|=1$. Passing to a subsequence, we have $\hat{w}_n \rightharpoonup \hat{w}$ in $\w$ and $\hat{w}_n \to \hat{w}$ in $L^p(\Omega)$, and $\hat{w}_n(x)\to 0$ for a.e.\ $x\in\Omega$. Since $\hat{w}_n^- \to 0$ in $\w$, we deduce that $\hat{w} \ge 0$ in $\Omega$. By \eqref{sps1} we have for all $n \in \N$, $v \in \w$
\beq \label{sps2}
\langle \fpl \hat{w}_n, v\rangle \le \int_\Omega \frac{f_1(x,w_n)}{\|w_n\|^{p-1}}v\,dx + \frac{\eps_n}{\|w_n\|^{p-1}} \|v\|.
\eeq
By ${\bf H}_1$ \ref{h11} and the definition of $f_1$, we have for all $n \in \N$
\[\int_\Omega \Big|\frac{f_1(x,w_n)}{\|w_n\|^{p-1}}\Big|^{p'}\,dx \le \int_\Omega \frac{C (1+ |w_n|^{p-1})^{p'}}{\|w_n\|^p}\,dx \le C.\]
By reflexivity of $L^{p'}(\Omega)$, passing to a subsequence we have
\beq \label{sps3}
\frac{f_1(\cdot, w_n)}{\|w_n\|^{p-1}} \rightharpoonup h_1 \text{ in } L^{p'}(\Omega).
\eeq 
Setting $v=\hat{w}_n-\hat{w} \in \w$ in \eqref{sps2}, we have for all $n\in\N$
\[\langle \fpl \hat{w}_n, \hat{w}_n-\hat{w} \rangle \le \int_\Omega \frac{f_1(x,w_n)}{\|w_n\|^{p-1}} (\hat{w}_n-\hat{w})\,dx + \frac{\eps_n}{\|w_n\|^{p-1}} \|\hat{w}_n-\hat{w}\|,\]
and the latter tends to $0$ as $n \to \infty$. By the $(S)_+$-property of $\fpl$, we have $\hat{w}_n \to \hat{w}$ in $\w$.
\vskip2pt
\noindent
We claim that there exists $\beta_1 \in L^{\infty}(\Omega)_+$ s.t.\ $\eta_1 \le \beta_1 \le \eta_2$ in $\Omega$ and $h_1=\beta_1 \hat{w}^{p-1}$. Indeed, set
\[\Omega^+=\big\{x \in \Omega:\,\hat{w}(x) >0\big\}, \quad \Omega^0=\big\{x \in \Omega:\,\hat{w}(x)=0\big\},\]
and for all $\eps>0$, $n \in \N$
\[\Omega_{\eps,n}^+=\Big\{x \in \Omega:\,w_n(x)>0, \ \eta_1(x)-\eps \le \frac{f_1(x,w_n)}{w_n^{p-1}(x)} \le \eta_2(x) + \eps \Big\}.\]
By ${\bf H}_1$ \ref{h12} we have $\chi_{\Omega_{\eps,n}^+} \to 1$ a.e.\ in $\Omega^+$ with dominated convergence, hence by \eqref{sps3} 
\[\chi_{\Omega_{\eps,n}^+} \frac{f_1(\cdot,w_n)}{\|w_n\|^{p-1}} \rightharpoonup h_1 \text{ in } L^{p'}(\Omega^+).\] 
So for all $\eps>0$, $n \in \N$ big enough, in $\Omega^+$ we have both $\hat{w}_n>0$ and
\[(\eta_1-\eps) \chi_{\Omega_{\eps,n}^+} \hat{w}_n^{p-1} \le \chi_{\Omega_{\eps,n}^+} \frac{f_1(\cdot,w_n)}{\|w_n\|^{p-1}} \le (\eta_2+\eps) \chi_{\Omega_{\eps,n}^+} \hat{w}_n^{p-1}.\] 
Passing to the limit as $n \to \infty$, we get in $\Omega^+$
\[(\eta_1-\eps) \hat{w}^{p-1} \le h_1 \le (\eta_2+\eps) \hat{w}^{p-1}.\]
Further, letting $\eps \to 0^+$, we obtain in $\Omega^+$
\[\eta_1 \hat{w}^{p-1} \le h_1 \le \eta_2 \hat{w}^{p-1}.\]
Similarly we see that $h_1=0$ in $\Omega^0$, and thus prove our claim.
\vskip2pt
\noindent
Passing to the limit in \eqref{sps2} as $n \to \infty$ we have weakly in $\Omega$
\beq \label{sps4}
\fpl \hat{w} = \beta_1(x) \hat{w}^{p-1}.
\eeq
By ${\bf H}_1$ \ref{h12} we have $\beta_1 \ge \lambda_1$ in $\Omega$, $\beta_1 \neq \lambda_1$. Hence, by Proposition \ref{fev} \ref{fev2} we have
\[\lambda_1(\beta_1) < \lambda_1(\lambda_1)=1.\]
So, by \eqref{sps4} $\hat{w} \in \w_+$ is a non-principal eigenfunction of $\fpl$ with weight $\beta_1 \in L^{\infty}(\Omega)_+$, against Proposition \ref{fev} \ref{fev1}. Such contradiction implies that $(w_n)$ is bounded in $\w$. By the bounded $(PS)$-property of $\Phi_1$, $(w_n)$ admits a convergent subsequence.
\vskip2pt
\noindent
The next step consists in proving that
\beq \label{sps5}
\lim_{\tau \to \infty} \Phi_1(\tau e_1)= - \infty,
\eeq
with $e_1\in{\rm int}(\cs_+)$ given by Proposition \ref{fev}. Indeed, by ${\bf H}_1$ \ref{h12} we have
\[\int_\Omega \eta_1(x)e_1^p\,dx > \lambda_1 \int_\Omega e_1^p\,dx=\lambda_1.\]
So we can fix $\eps>0$ s.t.\
\[\int_\Omega \eta_1(x)e_1^p\,dx > \lambda_1+\eps.\]
By ${\bf H}_1$ \ref{h12} again, there exists $T_{\eps}>\|u_0\|_{\infty}$ s.t.\ for a.e.\ $x\in\Omega$ and all 
$t > T_{\eps}$
\[f(x,t) \ge (\eta_1(x)-\eps)t^{p-1},\]
which, along with ${\bf H}_1$ \ref{h11}, implies
\begin{align*}
F_1(x,t) &= \int_0^{u_0} f(x,u_0)\,dx + \int_{u_0}^{T_{\eps}} f(x,\tau)\,d\tau + \int_{T_{\eps}}^t f(x,\tau)\,d\tau\\
&\ge f(x,u_0) u_0 - \int_{u_0}^{T_{\eps}} c_0 (1+ \tau^{p-1})\,d\tau + \int_{T_{\eps}}^t (\eta_1(x)-\eps)\tau^{p-1}\,d\tau\\
& \ge \frac{\eta_1(x)-\eps}{p} t^p - C_{\eps}, 
\end{align*}
for some $C_{\eps}>0$. Therefore we have for all $\tau >0$
\begin{align*}
\Phi_1(\tau e_1) &\le \frac{\tau^p \|e_1\|^p}{p} - \int_\Omega \Big(\frac{\eta_1(x)-\eps}{p} \tau^p e_1^p - C_{\eps}\Big)\,dx\\
&\le\frac{\lambda_1 \tau^p}{p} - \frac{\tau^p}{p} \int_\Omega \eta_1(x) e_1^p\,dx +  \frac{\eps \tau^p}{p} + C\\
&=\Big(\lambda_1- \int_\Omega \eta_1(x) e_1^p\,dx + \eps \Big)\frac{\tau^p}{p} + C,
\end{align*}
and the latter tends to $- \infty$ as $\tau \to \infty$, proving \eqref{sps5}.
\vskip2pt
\noindent
We define now a new truncation. Recalling that $0<u_0<a_+$ in $\Omega$, set for all $(x,t)\in\Omega\times\R$
\[\hat{f}_1(x,t)=
\begin{cases}
f(x,u_0) & \text{if $t<u_0(x)$,} \\
f(x,t) & \text{if $u_0(x)\le t \le a_+$,}\\
0  & \text{if $t > a_+$,}
\end{cases}\]
As usual, by ${\bf H}_1$ \ref{h11} \ref{h14} we see that $\hat{f}_1: \Omega \times \R \to \R$ satisfies ${\bf H}_0$. Setting
\[\hat{F}_1= \int_0^t \hat{f}_1(x,\tau)\,d\tau\]
and for all $u\in\w$
\[\hat{\Phi}_1(u)= \frac{\|u\|^p}{p} -  \int_{\Omega} \hat{F}_1(x,u)\,dx,\]
we deduce that $\hat{\Phi}_1 \in C^1(\w)$ is sequentially weakly lower semicontinuous. As in the proof of Lemma \ref{min} we see that $\hat{\Phi}_1$ is coercive in $\w$, so there exists $\hat{u}_1 \in \w$ s.t.\ 
\beq\label{sps6}
\hat{\Phi}_1(\hat{u}_1) = \inf_{u \in \w} \hat{\Phi}_1(u).
\eeq
In particular we have $\hat\Phi'_1(\hat{u}_1)=0$ in $W^{-s,p'}(\Omega)$, i.e., weakly in $\Omega$
\beq \label{sps7}
\fpl \hat{u}_1 = \hat{f}_1(x,\hat{u}_1).
\eeq
Testing \eqref{sps7} with $(u_0- \hat{u}_1)^+ \in \w$ and recalling that $u_0$ solves \eqref{dir},
\[\langle \fpl u_0 - \fpl \hat{u}_1, (u_0- \hat{u}_1)^+ \rangle = \int_\Omega \big(f(x,u_0)-\hat{f}_1(x,\hat{u}_1)\big) (u_0- \hat{u}_1)^+\,dx=0, \]
hence, by Proposition \ref{stm}, we have $\hat{u}_1 \ge u_0$ in $\Omega$. Similarly, testing \eqref{sps7} with $(\hat{u}_1-a_+)^+ \in \w$ we get $\hat{u}_1 \le a_+$ in $\Omega$. By construction, we can replace $\hat{f}_1$ with $f$ in \eqref{sps7}, so $\hat{u}_1$ is a solution of \eqref{dir}. By Proposition \ref{reg} we have $\hat u_1\in C_s^{\alpha}(\overline{\Omega})$. Also, since $\hat{u}_1 \ge u_0$, we easily get $\hat{u}_1 \in {\rm int}(\cs_+)$.
\vskip2pt
\noindent
To conclude, we distinguish two cases:
\begin{itemize}[leftmargin=1cm]
\item[$(a)$] If $\hat{u}_1 \neq u_0$, then by ${\bf H}_1$ \ref{h15} and Proposition \ref{comp} we have $\hat{u}_1> u_0$ in $\Omega$, and setting $u_1 = \hat{u}_1$ we conclude.
\item[$(b)$] If $\hat{u}_1 = u_0$, then set
\[\mathcal{U}=\big\{u \in \w \cap \cs: u<a_+ \ \text{in $\Omega$}\big\},\]
an open set in $\cs$ s.t.\ $u_0 \in \mathcal{U}$. By \eqref{sps6} we have for all $u \in \mathcal{U}$
\[\Phi_1(u)=\hat{\Phi}_1(u) \ge \hat{\Phi}_1(u_0)=\Phi_1(u_0).\]
So $u_0$ is a local minimizer of $\Phi_1$ in $\cs$, hence by Proposition \ref{svh} it is such as well in $\w$. Namely, there exists $\rho>0$ s.t.\ $\Phi_1(u)\ge\Phi_1(u_0)$ for all $u\in\overline{B}_\rho(u_0)$. Recalling \eqref{sps5}, we see that $\Phi_1$ exhibits the mountain pass geometry, while we have seen that it satisfies $(PS)$. By the mountain pass theorem (see for instance \cite[Theorem 5.40]{MMP}), there exists $u_1 \in \w \setminus \{u_0\}$ s.t.\ $\Phi'_1(u_1)=0$ in $W^{-s,p'}(\Omega)$ and $\Phi_1(u_1) \ge \Phi_1(u_0)$. So we have
\[\langle \fpl u_0 - \fpl u_1, (u_0- u_1)^+ \rangle = \int_\Omega (f(x,u_0)-f_1(x,u_1)) (u_0- u_1)^+\,dx=0, \]
hence, by Proposition \ref{stm}, we obtain $u_1 \ge u_0$ in $\Omega$. So $f_1(x,u_1)=f(x,u_1)$ for a.e.\ $x\in\Omega$, and we deduce that $u_1$ solves \eqref{dir}. Finally, by $u_1 \ge u_0$ and reasoning as above, we have $u_1 \in {\rm int}(\cs_+)$.
\end{itemize}
In both cases, the proof is concluded.
\end{proof}

\begin{remark}\label{strict}
One natural question about Lemma \ref{sps} is the following: under hyotheses ${\bf H}_1$, can we ensure that $u_1>u_0$ in $\Omega$? In general, the answer is negative, as in case $(b)$ above $u_1 > a_+$ may occur in $\Omega$, so we cannot use the quasi-monotonicity condition ${\bf H}_1$ \ref{h15} and Proposition \ref{comp}. Nevertheless, under the stricter assumption that
\[t\mapsto f(x,t)+c_2|t|^{p-2}t\]
is nondecreasing in $\R$ for a.e.\ $x\in\Omega$, we can prove that $u_1>u_0$ in $\Omega$.
\end{remark}

\noindent
Hypotheses ${\bf H}_1$ are somewhat symmetric, so we can argue as in Lemmas \ref{min}, \ref{sps} on the negative semiaxis and produce two negative solutions. As a whole, we have four constant sign solutions:

\begin{theorem}\label{css}
If ${\bf H}_1$ hold, then \eqref{dir} has at least four nontrivial solutions:
\begin{enumroman}
\item\label{css1} $u_0, u_1 \in {\rm int}(\cs_+)$, s.t.\ $0<u_0<a_+$, $u_0 \le u_1$ in $\Omega$,
\item\label{css2} $v_0, v_1 \in -{\rm int}(\cs_+)$, s.t.\ $a_-<v_0<0$, $v_1 \le v_0$ in $\Omega$.
\end{enumroman}
\end{theorem}

\noindent
Under the same hypotheses ${\bf H}_1$, we can reach more precise information about constant sign solutions (which will be useful in proving existence of a nodal solution in the next section):

\begin{lemma}\label{smp}
If ${\bf H}_1$ hold, then \eqref{dir} has a smallest positive solution $u_+ \in {\rm int}(\cs_+)$ and a biggest negative solution $v_-\in -{\rm int}(\cs_+)$.
\end{lemma}
\begin{proof}
First we set $\overline{u}=a_+ \in \widetilde{W}^{s,p}(\Omega)$ ($a_+>0$ introduced in ${\bf H}_1$ \ref{h14}), which is a supersolution of \eqref{dir}. Indeed, by ${\bf H}_1$ \ref{h14}, we have weakly in $\Omega$
\[\fpl\overline{u} = 0 = f(x,\overline{u}).\]
Now let $e_1 \in {\rm int}(\cs_+)$ be defined as in Proposition \ref{fev}, and $(\tau_n)$ be a decreasing sequence in $\R$ s.t.\ $\tau_n\to 0^+$, with $\tau_1>0$ small enough s.t.\ in $\Omega$
\[0 < \tau_1 e_1 \le \delta_0 < a_+, \quad \lambda_1(\tau_1 e_1)^{p-q} < c_1,\]
with $\delta_0,c_1>0$, $q\in(1,p)$ as in ${\bf H}_1$ \ref{h13} (and recalling from ${\bf H}_1$ \ref{h14} that $\delta_0<a_+$). For any $n\in\N$ set $\underline{u}_n=\tau_ne_1\in\w$. By ${\bf H}_1$ \ref{h13} and \eqref{evm} (with $m=1$) we have weakly in $\Omega$
\begin{align*}
\fpl \underline{u}_n &= \lambda_1 (\tau_n e_1)^{p-1} \\
&< c_1 (\tau_n e_1)^{q-1} \\
&\le f(x, \tau_n e_1) \\
&= f(x,\underline{u}_n).
\end{align*}
So, $\underline{u}_n$ is a (strict) subsolution of \eqref{dir}. Moreover, for all $n\in\N$ we have $\underline{u}_n\le\overline{u}$ in $\Omega$, so $(\underline{u}_n,\overline{u})$ is a sub-supersolution pair. By Proposition \ref{ext}, for all $n\in\N$ we can find a function $\check{u}_n \in \w_+$ s.t.\ 
\beq\label{smp1}
\check{u}_n = \min \mathcal{S}(\underline{u}_n,\overline{u}).
\eeq
In particular, for all $n\in\N$ we have weakly in $\Omega$
\beq\label{smp2}
\fpl \check{u}_n = f(x,\check{u}_n).  
\eeq
Testing \eqref{smp2} with $\check{u}_n\in\w$ and using ${\bf H}_1$ \ref{h11}, we get for all $n\in\N$
\begin{align*}
\|\check{u}_n \|^p &= \langle \fpl \check{u}_n, \check{u}_n \rangle \\
&= \int_\Omega  f(x,\check{u}_n) \check{u}_n\,dx\\
&\le \int_\Omega c_0 (\check{u}_n + \check{u}_n^p)\,dx\\
&\le c_0 (a_+ + a_+^p) |\Omega|.
\end{align*}
Hence $(\check{u}_n)$ is bounded in $\w$. Passing to a subsequence, $\check{u}_n \rightharpoonup u_+$ in $\w$, $\check{u}_n \to u_+$ in both $L^p(\Omega)$ and $L^1(\Omega)$, for some $u_+ \in \w_+$. Testing \eqref{smp2} with $(\check{u}_n - u_+)$ leads to
\begin{align*}
\langle \fpl \check{u}_n, \check{u}_n - u_+ \rangle &=  \int_\Omega f(x,\check{u}_n) (\check{u}_n- u_+)\,dx\\ 
&\le \int_\Omega c_0 (1+ \check{u}_n^{p-1}) (\check{u}_n- u_+)\,dx\\
&\le c_0 \big(\|\check{u}_n- u_+\|_1 +  \|\check{u}_n\|_p^{p-1} \|\check{u}_n- u_+\|_p\big) 
\end{align*}
(where we have used ${\bf H}_1$ \ref{h11} and H\"{o}lder's inequality), and the latter tends to $0$ as $n \to \infty$. By the $(S)_+$-property of $\fpl$, we get $\check{u}_n \to u_+$ in $\w$. So we may pass to the limit as $n \to \infty$ in \eqref{smp2}, and we have weakly in $\Omega$
\[\fpl u_+ = f(x,u_+). \]
Therefore $u_+ \in  C_s^{\alpha}(\overline{\Omega})_+$ is a solution of \eqref{dir} (see Proposition \ref{reg}). The next step consists in proving that
\beq \label{smp3}
u_+ \neq 0.
\eeq
We introduce the auxiliary problem (with $c_1, q$ as in ${\bf H}_1$ \ref{h13}):
\beq \label{smp4}
\begin{cases}
\fpl w = c_1 w^{q-1} & \text{in $\Omega$,} \\
w\ge 0 & \text{in $\Omega$,}\\
w=0 & \text{in $\Omega^c$.}
\end{cases}
\eeq
The energy functional $\Psi \in C^1(\w)$ of \eqref{smp4} is defined by setting for all $u\in\w$
\[\Psi(u)=\frac{\|u\|^p}{p}-c_1\frac{\|u^+\|_q^q}{q},\]
and it is coercive and sequentially weakly lower semicontinuous. Besides, for all $\tau>0$ we have
\[\Psi(\tau e_1) = \frac{\lambda_1}{p}\tau^p-c_1\frac{\|e_1\|_q^q}{q}\tau^q,\]
and the latter is negative for $\tau>0$ small enough. Hence there exists $w \in \w$ s.t.\ 
\[\Psi(w)=\inf_{u \in \w} \Psi(u)<0.\]
In particular, we have $\Psi'(w)=0$ in $W^{-s,p'}(\Omega)$. Testing such relation with $-w^-\in\w$ and using \eqref{pnp}, we get
\begin{align*}
\|w^-\|^p &\le \langle\fpl w,-w^-\rangle \\
&= \int_\Omega (w^+)^{q-1}(-w^-)\,dx = 0,
\end{align*}
hence $w\in\w_+\setminus\{0\}$ and it solves \eqref{smp4}. By Propositions \ref{reg} and \ref{max} (with $g(t)=0$) we easily get $w \in {\rm int}(\cs_+)$. Let $\theta \in (0,1)$ be s.t.\ $0 < \theta w \le \delta_0$ in $\Omega$, then set 
$\check{w}=\theta w \in {\rm int}(\cs_+)$. For all $n \in \N$ we have in $\Omega$
\beq \label{smp5}
\check{w} \le \check{u}_n.
\eeq
Indeed, since $\check{w}, \check{u}_n \in {\rm int}(\cs_+)$, for all $\sigma>0$ small enough we have $\check{u}_n-\sigma\check{w}\in{\rm int}(\cs_+)$, hence in $\Omega$
\[\sigma \check{w} \le \check{u}_n.\]
Arguing by contradiction, let $\sigma\in(0,1)$ be maximal in the inequality above. In particular we have $\sigma \check{w} \neq \check{u}_n$ and by \eqref{smp4} we have weakly in $\Omega$
\begin{align*}
\fpl (\sigma \check{w}) + c_2 (\sigma \check{w})^{p-1} &= (\sigma \theta)^{p-1} c_1 w^{q-1} + c_2 (\sigma \check{w})^{p-1}\\ 
&< c_1 (\sigma \check{w})^{q-1} + c_2 (\sigma \check{w})^{p-1}\\
& \le f(x,\sigma \check{w}) + c_2 (\sigma \check{w})^{p-1}\\ 
&\le f(x,\check{u}_n) + c_2 \check{u}_n^{p-1}\\ 
&= \fpl \check{u}_n + c_2 \check{u}_n^{p-1},
\end{align*}
where we have used ${\bf H}_1$ \ref{h13} \ref{h15} and \eqref{smp2}. We apply Proposition \ref{comp} (with $g(t)=c_2|t|^{p-2}t$) and find that $\check{u}_n - \sigma \check{w} \in {\rm int}(\cs_+)$, hence as above we can find $\sigma' \in (\sigma,1)$ s.t.\ in $\Omega$
\[\sigma' \check{w} < \check{u}_n,\]
against maximality. Thus we have \eqref{smp5}. Passing to the limit in \eqref{smp5} as $n \to \infty$ and using \eqref{smp1} we have $\check{w} \le u_+$ in $\Omega$, which proves \eqref{smp3}.
\vskip2pt
\noindent
Recall that $0 \le u_+ \le a_+$ in $\Omega$, hence by ${\bf H}_1$ \ref{h15} we have
\[\begin{cases}
\fpl u_+ + c_2(u_+)^{p-1} = f(x,u_+)+c_2(u_+)^{p-1} \ge 0 & \text{weakly in $\Omega$} \\
u_+ \ge 0 & \text{in $\R^N$}.
\end{cases}\]
By Proposition \ref{max} (with $g(t)=c_2|t|^{p-2}t$) we have $u_+ \in {\rm int}(\cs_+)$.
\vskip2pt
\noindent
Finally, we need to show that $u_+$ is the smallest solution of \eqref{dir}. Indeed, let $u \in {\rm int}(\cs_+)$ be any positive solution of \eqref{dir}. By Lemma \ref{zero} we have
\[\check{u} = \min\{u,\overline{u}\} \in \w,\]
and by Proposition \ref{sss} $\check{u}$ is a supersolution of \eqref{dir}. Also, for all $n\in\N$ big enough we have $\underline{u}_n\le\check{u}\le\overline{u}$ in $\Omega$. By Proposition \ref{ext} there exists $\hat{u}_n \in \mathcal{S}(\underline{u}_n, \check{u})$, i.e., $\hat{u}_n \in {\rm int}(\cs_+)$ solves \eqref{dir} and in $\Omega$ we have
\[0 < \underline{u}_n \le \hat{u}_n \le \check{u} \le \overline{u}.\]
By \eqref{smp1} we have in $\Omega$
\[\check{u}_n \le \hat{u}_n \le u.\]
Passing to the limit as $n \to \infty$ we have $u_+ \le u$ in $\Omega$, which proves our claim.
\vskip2pt
\noindent
Reasoning similarly on the negative semiaxis, we produce a biggest negative solution $v_- \in -{\rm int}(\cs_+)$.
\end{proof}

\section{Nodal solution}\label{sec4}

\noindent
In this section we seek a nodal (i.e., sign-changing) solution of \eqref{dir}, taking values between the extremal constant sign solutions seen in Lemma \ref{smp}. We are going to exploit some Morse theory (computation of critical groups of the energy functional). For such a purpose, we need to stregthen slightly our hypotheses on the reaction $f$, adding a reverse Ambrosetti-Rabinowitz condition near the origin:

\begin{itemize}[leftmargin=1cm]
\item[${\bf H}_2$] $f:\Omega\times\R\to\R$ is a Carath\'{e}odory mapping satisfying \ref{h11}-\ref{h15} as in ${\bf H}_1$ and
\begin{enumroman}
\setcounter{enumi}{5}
\item\label{h16} there exists $\mu\in(1,p)$ s.t.\ for a.e.\ $x\in\Omega$ and all $|t|\le\delta_0$
\[\mu F(x,t) \ge f(x,t)t.\]
\end{enumroman}
\end{itemize}

\noindent
Clearly, ${\bf H}_2$ implies ${\bf H}_1$ (and hence ${\bf H}_0$). Thus, all results of Sections \ref{sec2}, \ref{sec3} apply.

\begin{example}\label{exnew}
The function $f\in C^0(\R)$ introduced in Example \ref{exa} satisfies ${\bf H}_2$ as well. Indeed, a straightforward calculation leads to
\[f(t) = t^\frac{1}{2}+{\bf o}(t^\frac{1}{2}), \ F(t) = \frac{2}{3}t^\frac{3}{2}+{\bf o}(t^\frac{3}{2}),\]
as $t\to 0^+$. Hence, taking $\mu\in(3/2,2)$ (recall that $p=2$), for all $t>0$ small enough we have
\[\mu F(t) \ge f(t)t\]
(the case $t<0$ is studied similarly). Thus $f$ satisfies ${\bf H}_2$ \ref{h16}.
\end{example}

\noindent
By Lemma \ref{smp}, under ${\bf H}_2$ problem \eqref{dir} admits the extremal constant sign solutions $u_+ \in {\rm int}(\cs_+)$, $v_- \in -{\rm int}(\cs_+)$. Set for all $(x,t) \in \Omega \times \R$
\[\tilde{f}(x,t)=
\begin{cases}
f(x,v_-) & \text{if $t<v_-$,} \\
f(x,t) & \text{if $v_- \le t \le u_+$,}\\
f(x,u_+) & \text{if $t>u_+$,}
\end{cases}\]
and accordingly set
\[\tilde{F}(x,t) = \int_0^t \tilde{f}(x,\tau)\,d\tau.\]
Further, set for all $u\in\w$
\[\tilde{\Phi}(u)= \frac{\|u\|^p}{p} - \int_{\Omega} \tilde{F}(x,u)\,dx.\]
By ${\bf H}_2$ \ref{h11} we see that $\tilde{f}:\Omega\times\R\to\R$ satisfies ${\bf H}_0$, hence $\tilde{\Phi} \in C^1(\w)$ with derivative given for all $u,v\in\w$ by
\[\langle \tilde{\Phi}'(u),v \rangle = \langle \fpl u,v \rangle - \int_\Omega \tilde{f}(x,u)v\,dx.\]
By ${\bf H}_2$ \ref{h14} we have $\tilde{\Phi}'(0)=0$ in $W^{-s,p'}(\Omega)$. Consistently with the general assumptions that \eqref{dir} has finitely many solutions, without loss of generality we may assume that $0$ is an isolated critical point of $\tilde{\Phi}$ (see Section \ref{sec2}).
\vskip2pt
\noindent
In the next lemma we compute the critical groups of $\tilde{\Phi}$ at $0$: 

\begin{lemma}\label{cgz}
If ${\bf H}_2$ hold, then $C_k(\tilde{\Phi},0)=0$ for all $k \in \N$.
\end{lemma}
\begin{proof}
Fix any $r \in (p,p_s^*)$. Combining ${\bf H}_2$ \ref{h11} \ref{h13}, we can find $C_1,C_2>0$ s.t.\ for a.e.\ $x\in\Omega$ and all $t \in \R$
\beq \label{cgz1}
\tilde{F}(x,t) \ge C_1 |t|^q - C_2 |t|^r.
\eeq
Now fix $u \in \w \setminus \{0\}$. By \eqref{cgz1} we have for all $\tau>0$
\begin{align*}
\tilde{\Phi}(\tau u) &\le \frac{\tau^p \|u\|^p}{p} - \int_\Omega \big(C_1 \tau^q |u|^q - C_2 \tau^r |u|^r\big)\,dx \\
&= \frac{\tau^p \|u\|^p}{p} - C_1 \tau^q \|u\|_q^q + C_2 \tau^r \|u\|_r^r,
\end{align*}
and the latter is negative for $\tau>0$ small enough (since $q<p<r$). By Lemma \ref{zero} we have
\[w_0=\min\{u_+,-v_-, \delta_0\} \in \w.\]
We begin studying the behavior of $\tilde\Phi$ near $0$. First, we claim that there exists $\rho>0$ s.t.\ for all $u \in \overline{B}_{\rho}(0) \setminus \{0\}$ with $\tilde{\Phi}(u)=0$ we have
\beq \label{cgz2}
\restr{\frac{d}{d\tau}\tilde{\Phi}(\tau u)}{\tau=1} >0.
\eeq
Indeed, fix $u \in \w\setminus \{0\}$ with $\tilde{\Phi}(u)=0$. Using ${\bf H}_2$ \ref{h16} and the definition of $\tilde f$, we compute for all $\tau>0$
\begin{align}\label{cgz3}
\restr{\frac{d}{d\tau}\tilde{\Phi}(\tau u)}{\tau=1}&= \langle \tilde{\Phi}'(u),u\rangle \\
\nonumber &= \|u\|^p - \int_\Omega \tilde{f}(x,u)u\,dx \\
\nonumber &= \Big(1-\frac{\mu}{p}\Big) \|u\|^p + \int_\Omega \big(\mu \tilde{F}(x,u)-\tilde{f}(x,u)u\big)\,dx \\ 
\nonumber &\ge \Big(1-\frac{\mu}{p}\Big) \|u\|^p + \int_{\{|u| \ge w_0\}} \big(\mu \tilde{F}(x,u)-\tilde{f}(x,u)u\big)\,dx.
\end{align}
There exists $C>0$ (independent of $u$) s.t.\
\beq\label{cgz4}
\int_{\{|u| \ge w_0\}} \big(\mu \tilde{F}(x,u)-\tilde{f}(x,u)u\big)\,dx \ge - C \|u\|_r^r.
\eeq
Indeed, fix $x \in \Omega$ s.t.\ $u(x) \ge w_0(x)$. Three cases may occur:
\begin{itemize}[leftmargin=1cm]
\item[$(a)$] If $\delta_0 \le u(x) \le u_+(x)$, then by ${\bf H}_2$ \ref{h11} and \eqref{cgz1} we have
\begin{align*}
\mu \tilde{F}(x,u)-\tilde{f}(x,u)u &= \mu F(x,u)-f(x,u)u\\
& \ge \mu C (u^q-u^r)-C(u+u^p)\\
& \ge -C(u+u^p+u^r) \\
&\ge -C \Big(\frac{u^r}{\delta_0^{r-1}} + \frac{u^r}{\delta_0^{r-p}} + u^r \Big)\\ 
&\ge -C u^r.
\end{align*}
\item[$(b)$]  If $u_+(x)\le u(x) \le \delta_0 $, then, recalling the definition of $\tilde{f}$, that $\mu >1$, and that $f(x,u_+)\ge 0$ (by ${\bf H}_2$ \ref{h13}), we obtain
\begin{align*}
\mu \tilde{F}(x,u)-\tilde{f}(x,u)u &= \mu \int_0^{u_+} f(x,t)\,dt + \mu \int_{u_+}^u f(x,u_+)\,dt - f(x,u_+)u\\
&= \mu F(x,u_+) + \mu f(x,u_+) (u-u_+) - f(x,u_+)u\\
& \ge \mu F(x,u_+) - f(x,u_+)u_+ \ge 0.
\end{align*}
\item[$(c)$] If $-v_-(x)\le u(x) \le \min\{u_+(x),\delta_0\}$, simply use ${\bf H}_2$ \ref{h16}.
\end{itemize}
Similarly we deal with the case $u(x)<0$. Therefore, plugging \eqref{cgz4} into \eqref{cgz3} we get
\[\restr{\frac{d}{d\tau}\tilde{\Phi}(\tau u)}{\tau=1} \ge \Big(1-\frac{\mu}{p}\Big) \|u\|^p - C \|u\|^r,\]
and the latter is positive whenever $\|u\|>0$ is small enough. So we have \eqref{cgz2}.
\vskip2pt
\noindent
Our next claim is that for all $u\in\overline{B}_\rho(0)\setminus\{0\}$ s.t.\ $\tilde\Phi(u)=0$ and all $\tau\in[0,1]$
\beq \label{cgz5}
\tilde{\Phi}(\tau u) \le 0.
\eeq
Arguing by contradiction, let $\tau_1 \in (0,1)$ s.t.\ $\tilde{\Phi}(\tau_1 u) >0$. By the mean value theorem there exists $\tau_2 \in (\tau_1,1]$ minimal s.t.\ $\tilde{\Phi}(\tau_2 u)=0$. So $\tilde{\Phi}(\tau u) >0$ for all $\tau \in [\tau_1, \tau_2)$. Set $w=\tau_2 u$, then $w \in \overline{B}_{\rho}(0) \setminus \{0\}$ and $\tilde{\Phi}(w)=0$. By \eqref{cgz2} we have
\[\restr{\frac{d}{d\tau}\tilde{\Phi}(\tau w)}{\tau=1} > 0.\]
Besides, since $\tilde{\Phi}(\tau w) >0$ for all $\tau \in(\tau_1/\tau_2,1)$ we have
\[\restr{\frac{d}{d\tau}\tilde{\Phi}(\tau w)}{\tau=1} \le 0,\]
a contradiction. So \eqref{cgz5} is proved.
\vskip2pt
\noindent
Taking $\rho>0$ even smaller if necessary, we have
\[K(\tilde\Phi)\cap\overline{B}_\rho(0) = \{0\}.\]
Set
\[A = \big\{u \in \overline{B}_{\rho}(0): \tilde{\Phi}(u) \le 0\big\}.\]
Clearly $0 \in A$. Plus, $A$ is a star-shaped set centered at $0$. Indeed, for all $u \in A \setminus \{0\}$, $\tau \in [0,1]$ we have $\tilde{\Phi}(\tau u) \le 0$. Otherwise, there would exist $0<\tau_1<\tau_2<1$ s.t.\ 
\[\tilde{\Phi}(\tau_1 u)>0=\tilde{\Phi}(\tau_2 u),\]
against \eqref{cgz5}. By \cite[Remark 6.23]{MMP}, the set $A$ is contractible. Now consider $u \in \overline{B}_{\rho}(0) \setminus A$, i.e., satisfying $\tilde{\Phi}(u)>0$. As seen above, $\tilde{\Phi}(\tau u)<0$ for all $\tau \in (0,1)$ small enough, so there exists $\tau \in (0,1)$ s.t.\
\[\tilde{\Phi}(\tau u)=0.\]
We claim that such $\tau \in (0,1)$ is unique. Arguing by contradiction, let $0<\tau_1<\tau_2<1$ be s.t.\
\[\tilde{\Phi}(\tau_1 u) = \tilde\Phi(\tau_2 u) = 0.\]
By \eqref{cgz5} we have for all 
$\sigma \in [0,1]$
\[\tilde{\Phi}(\sigma \tau_2 u) \le 0.\]
Now set for all $\sigma\in[0,1]$
\[g(\sigma) = \tilde{\Phi}(\sigma \tau_2 u),\]
so that the map $g\in C^1([0,1])$ attains its maximum at $\sigma=\tau_1/\tau_2 \in (0,1)$, hence
\[g'\Big(\frac{\tau_1}{\tau_2}\Big) = 0.\]
So we have
\[\restr{\frac{d}{d\tau}\tilde{\Phi}(\tau u)}{\tau=\tau_1} = \frac{1}{\tau_2}g'\Big(\frac{\tau_1}{\tau_2}\Big) = 0,\]
against \eqref{cgz2}. By the implicit function theorem (see \cite[Theorem 7.3]{MMP}), we can construct a continuous map $\hat{\tau}: \overline{B}_{\rho}(0) \setminus A \to (0,1)$ s.t.\ for all $u \in \overline{B}_{\rho}(0) \setminus A$, $\tau \in (0,1)$
\[\tilde{\Phi}(\tau u) =
\begin{cases}
<0 & \text{if $0<\tau<\hat{\tau}(u)$} \\
=0 & \text{if $\tau=\hat{\tau}(u)$} \\
>0 & \text{if $\hat{\tau}(u)<\tau<1$.}
\end{cases}\]
Set for all $u \in \overline{B}_{\rho}(0) \setminus \{0\}$
\[\hat{h}(u)=
\begin{cases}
u & \text{if $u \in A \setminus \{0\}$} \\
\hat{\tau}(u) u & \text{if $u \in  \overline{B}_{\rho}(0) \setminus A$.}
\end{cases}\]
Then $\hat{h}: \overline{B}_{\rho}(0) \setminus \{0\} \to A \setminus \{0\}$ is a continuous retraction. Since ${\rm dim}(\w)=\infty$, the set $\overline{B}_{\rho}(0) \setminus \{0\}$ is contractible \cite[Example 6.45 $(b)$]{MMP}. Then, being a retract of $\overline{B}_{\rho}(0) \setminus \{0\}$, $A \setminus \{0\}$ is contractible as well. So, by \cite[Propositions 6.24, 6.25]{MMP} we have for all $k\ge 0$
\[C_k(\tilde{\Phi},0)=H_k(A, A \setminus \{0\})=0,\]
which concludes the proof.
\end{proof}

\noindent
Finally, we can prove our complete multiplicity result:

\begin{theorem}\label{nod}
If ${\bf H}_2$ hold, then \eqref{dir} has at least five nontrivial solutions:
\begin{enumroman}
\item\label{nod1} $u_0, u_1 \in {\rm int}(\cs_+)$, s.t.\ $0<u_0<a_+$, $u_0 \le u_1$ in $\Omega$,
\item\label{nod2} $v_0, v_1 \in -{\rm int}(\cs_+)$, s.t.\ $a_-<v_0<0$, $v_1 \le v_0$ in $\Omega$,
\item\label{nod3} $\tilde{u} \in \cs \setminus \{0\}$ nodal, s.t.\ $v_0 \le \tilde{u} \le u_0$ in $\Omega$.
\end{enumroman}
\end{theorem}
\begin{proof}
From Theorem \ref{css} we have \ref{nod1}, \ref{nod2}. So, there remains to prove \ref{nod3}.
\vskip2pt
\noindent
By Lemma \ref{smp} we know that \eqref{dir} admits extremal constant sign solutions $u_+\in{\rm int}(\cs_+)$, $v_-\in{\rm int}(\cs_+)$. Without loss of generality and consistently with \ref{nod1}, \ref{nod2}, we may assume that
\[u_+ = u_0, \quad v_- = v_0.\]
In particular, that implies that no nontrivial, constant sign solution may exist in the set ${\mathcal S}(v_0,u_0)$. Now set for all $(x,t) \in \Omega \times \R$
\[\tilde{f}_{\pm}(x,t)=\tilde{f}(x, \pm t^{\pm}),\]
and
\[\tilde{F}_{\pm}= \int_0^t \tilde{f}_{\pm}(x,\tau)\,d\tau.\]
Further, set for all $u\in\w$
\[\tilde{\Phi}_{\pm}(u)= \frac{\|u\|^p}{p}-\int_{\Omega} \tilde{F}_{\pm}(x,u)\,dx.\]
Clearly $\tilde f_\pm:\Omega\times\R\to\R$ satisfies ${\bf H}_0$, so $\tilde{\Phi}_{\pm}\in C^1(\w)$ with derivative given for all $u,v\in\w$ by
\[\langle\tilde\Phi_\pm(u),v\rangle = \langle\fpl u,v\rangle-\int_\Omega \tilde f_\pm(x,u)v\,dx.\]
We now focus on $\tilde\Phi_+\in C^1(\w)$, proving that for all $u\in K(\tilde\Phi_+)$ we have in $\Omega$
\beq\label{nod4}
0 \le u \le u_0.
\eeq
Indeed, we have $\tilde{\Phi}'_+(u)=0$ in $W^{-s,p'}(\Omega)$, which rephrases as
\beq \label{nod5}
\fpl u = \tilde{f}_+(x,u)
\eeq
weakly in $\Omega$. Testing \eqref{nod5} with $(u-u_0)^+ \in \w$ and recalling that $u_0$ solves \eqref{dir}, we get
\[\langle \fpl u - \fpl u_0, (u- u_0)^+ \rangle =  \int_\Omega \big(\tilde{f}_+(x,u)-f(x,u_0)\big) (u- u_0)^+\,dx=0.\]
By Proposition \ref{stm}, we have $u \le u_0$ in $\Omega$. In a similar way, testing \eqref{nod5} with $-u^- \in \w$ and using \eqref{pnp} we have
\[\|u^-\|^p \le \langle\fpl u,-u^-\rangle = \int_\Omega\tilde f_+(x,u)(-u^-)\,dx = 0,\]
hence $u\ge 0$ in $\Omega$. So \eqref{nod4} is proved. More precisely, we have
\beq \label{nod6}
K(\tilde{\Phi}_+)=\{0,u_0\}.
\eeq
Indeed, by construction of $\tilde f_+$, clearly we have in $W^{-s,p'}(\Omega)$
\[\tilde\Phi'_+(0) = \tilde\Phi'_+(u_0) = 0.\]
Vice versa, let $u\in K(\tilde{\Phi}_+)$, i.e., $u$ satisfies \eqref{nod5}. By \eqref{nod4} we can replace $\tilde f_+$ with $f$ in \eqref{nod5}, hence $u\in C^\alpha_s(\overline\Omega)$ solves \eqref{dir} (see Proposition \ref{reg}). Assume $u \neq 0$. Then by ${\bf H}_2$ \ref{h15} we have
\[\begin{cases}
\fpl u + c_2 u^{p-1} = f(x,u)+c_2 u^{p-1} \ge 0 & \text{in $\Omega$} \\
u \ge 0 & \text{in $\Omega$.}
\end{cases}\]
By Proposition \ref{max} (with $g(t)=c_2|t|^{p-2}t$), we get $u \in {\rm int}(\cs_+)$. Therefore $u$ is a positive solution of \eqref{dir} s.t.\ $u \le u_0$. By minimality we deduce $u=u_0$, which proves \eqref{nod6}.
\vskip2pt
\noindent
Next, we claim that $u_0$ is a local minimizer of $\tilde\Phi$. Indeed, since $\tilde f_+:\Omega\times\R\to\R$ is bounded, then $\tilde\Phi_+$ is coercive and sequentially weakly lower semicontinuous. So there exists $\tilde{u}_+ \in \w$ s.t.\ 
\beq\label{nod7}
\tilde{\Phi}_+ (\tilde{u}_+)= \inf_{u \in \w} \tilde{\Phi}_+(u) = \tilde m_+.
\eeq
Let once again $e_1 \in {\rm int}(\cs_+)$ be as in Proposition \ref{fev}. For all $\tau>0$ small enough we have in $\Omega$
\[0 < \tau e_1 \le \min\{u_0, \delta_0\}.\]
Therefore, by ${\bf H}_2$ \ref{h13} we have in $\Omega$
\[\tilde{f}_+(x, \tau e_1)= f(x, \tau e_1) \ge c_1(\tau e_1)^{q-1}.\]
Hence, for all $\tau >0$ small enough
\begin{align*}
\tilde{\Phi}_+(\tau e_1) &\le \frac{\tau^p \|e_1\|^p}{p} - \int_\Omega \frac{c_1}{q} (\tau e_1)^q\,dx\\ 
&= \frac{\lambda_1}{p} \tau^p - \frac{c_1}{q} \tau^q  \|e_1\|_q^q,
\end{align*}
and the latter is negative for all $\tau >0$ small enough (since $q<p$). So, in \eqref{nod7} we have $\tilde m_+<0$, hence $\tilde u_+\neq 0$. Therefore we have $\tilde{u}_+ \in K(\tilde{\Phi}_+) \setminus \{0\}$, which by \eqref{nod6} implies $\tilde{u}_+ = u_0$. Further, since $u_0 \in {\rm int}(\cs_+)$, we can find $\rho>0$ s.t.\ for all $u \in \w \cap \cs$ with $\|u-u_0\|_{0,s} \le \rho$ we have $u>0$ in $\Omega$. So for any such $u$ we get
\[\tilde{\Phi}(u)=\tilde{\Phi}_+(u) \ge \tilde{\Phi}_+(u_0)= \tilde{\Phi}(u_0),\]
i.e., $u_0$ is a local minimizer of $\tilde{\Phi}$ in $\cs$. By Proposition \ref{svh}, $u_0$ is as well a local minimizer of $\tilde\Phi$ in $\w$.
\vskip2pt
\noindent
Reasoning in a similar way, we see that
\[K(\tilde{\Phi}_-)=\{0,v_0\},\]
and that $v_0$ is a local minimizer of $\tilde\Phi$.
\vskip2pt
\noindent
Now let us turn to the functional $\tilde\Phi\in C^1(\w)$. This is coercive, hence it satisfies $(PS)$ (see Section \ref{sec2}) and admits two distinct local minimizers $u_0,v_0\in\w$. Without loss of generality, we may assume that $K(\tilde\Phi)$ is a finite set, hence in particular that $u_0$ and $v_0$ are strict local minimizers. As in the proof of \eqref{nod4}, we see that for all $u\in K(\tilde\Phi)$ we have in $\Omega$
\beq\label{nod8}
v_0 \le u \le u_0.
\eeq
Now set
\[\Gamma = \big\{\gamma\in C([0,1],\w):\,\gamma(0)=u_0,\,\gamma(1)=v_0\big\}.\]
Then we have
\[c = \inf_{\gamma\in\Gamma}\max_{t\in[0,1]}\tilde\Phi(\gamma(t)) > \max\{\tilde\Phi(u_0),\,\tilde\Phi(v_0)\big\}.\]
By Hofer's version of the mountain pass theorem and the characterization of critical groups at mountain pass type critical points (see \cite[Theorem 6.99, Proposition 6.100]{MMP}), there exists $\tilde u\in K_c(\tilde\Phi)$ s.t.\
\beq \label{nod9}
C_1(\tilde{\Phi},\tilde{u}) \neq 0. 
\eeq
Then we have weakly in $\Omega$
\[\fpl\tilde u = \tilde f(x,\tilde u).\]
By \eqref{nod8} we have in $\Omega$
\[v_0 \le \tilde u \le u_0.\]
Therefore, by construction of $\tilde f$, we see that $\tilde u$ solves \eqref{dir}. Hence, by Proposition \ref{reg}, we have $\tilde u\in C^\alpha_s(\overline\Omega)$.
\vskip2pt
\noindent
There remains to prove that $\tilde u$ is nodal. Recall that $u_0,v_0\in K(\tilde\Phi)$ are strict local minimizers and isolated critical points. So, by \cite[Example 6.45 (a)]{MMP}, we obtain 
\[C_1(\tilde{\Phi},u_0) = C_1(\tilde{\Phi},v_0)=0. \]
Besides, by Lemma \ref{cgz} we have
\[C_1(\tilde{\Phi},0)=0.\]
So, \eqref{nod9} implies
\[\tilde u\in K(\tilde\Phi)\setminus\{0,u_0,v_0\}.\]
Assume now that $\tilde u\ge 0$ in $\Omega$. Then, by Proposition \ref{max} we would have $\tilde u\in{\rm int}(\cs_+)$ and $\tilde u\le u_0$, against minimality of $u_0$. Similarly, assuming that $\tilde u\le 0$ in $\Omega$ leads to a contradiction. Thus, $\tilde u$ changes sign in $\Omega$, which completes the argument for \ref{nod3}.
\end{proof}

\begin{remark}
As seen in Remark \ref{strict}, in general we cannot be sure that all solutions lie in the order interval $[a_-,a_+]$, so we cannot use Proposition \ref{comp} and get strict inequalities. These can be retrieved by strengthening the quasi-monotonicity hypothesis ${\bf H}_2$ \ref{h15}.
\end{remark}

\vskip4pt
\noindent
{\bf Acknowledgement.} Both authors are members of GNAMPA (Gruppo Nazionale per l'Analisi Matematica, la Probabilit\`a e le loro Applicazioni) of INdAM (Istituto Nazionale di Alta Matematica 'Francesco Severi') and are supported by the research project \emph{Evolutive and Stationary Partial Differential Equations with a Focus on Biomathematics}, funded by Fondazione di Sardegna (2019). A.\ Iannizzotto is also supported by the grant PRIN-2017AYM8XW: \emph{Nonlinear Differential Problems via Variational, Topological and Set-valued Methods}.

\end{document}